\documentclass[reqno,11pt]{article} 
\usepackage{amsmath,amssymb,amsthm,amscd,verbatim}
\usepackage{stmaryrd}
\usepackage{graphicx,caption,subcaption}
\usepackage{color}
\usepackage{breqn}
\usepackage{url}
\usepackage{latexsym}
\usepackage{amsmath}
\usepackage{amsfonts}
\usepackage{amssymb}
\usepackage{amsthm}
\usepackage{mathrsfs}
\usepackage{enumerate}

\usepackage{enumerate}

\usepackage{amsmath}
\usepackage{amstext,amssymb,amsfonts}
\usepackage{fullpage}
\usepackage{bm}

\usepackage{todonotes}

\newcounter{mynotes}
\usepackage{nameref}
\usepackage[pagebackref,colorlinks,linkcolor=blue,filecolor = blue,
citecolor = blue, urlcolor = blue,pdfstartview=FitH]{hyperref}
\usepackage[capitalize]{cleveref}

\usepackage{amsthm}
\usepackage{thmtools}

\declaretheorem[within=section]{theorem}

\declaretheorem[sibling=theorem]{lemma}

\declaretheorem[sibling=theorem]{definition}
\declaretheorem[sibling=theorem]{Lemma+Definition}

\declaretheorem[sibling=theorem]{proposition}

\crefname{proposition}{Proposition}{Propositions}
\crefname{conjecture}{Conjecture}{Conjectures}
\crefname{claim}{Claim}{Claims}
\crefname{remark}{Remark}{Remarks}

\newcounter{termcounter}
\renewcommand{\thetermcounter}{\Alph{termcounter}}
\crefname{term}{term}{terms}
\creflabelformat{term}{\textup{(#2#1#3)}}

\makeatletter
\def\term{\@ifnextchar[\term@optarg\term@noarg}
\def\term@optarg[#1]#2{%
  \textup{(#1)}%
  \def\@currentlabel{#1}%
  \def\cref@currentlabel{[][2147483647][]#1}%
  \cref@label[term]{#2}}
\def\term@noarg#1{%
  \refstepcounter{termcounter}%
  \textup{(\thetermcounter)}%
  \cref@label[term]{#1}}
\makeatother

\newcommand{\ignore}[1]{}

\definecolor{DSred}{rgb}{1,0,0}

\renewcommand{\leq}{\leqslant}
\renewcommand{\geq}{\geqslant}
\renewcommand{\ge}{\geqslant}
\renewcommand{\le}{\leqslant}
\renewcommand{\epsilon}{\varepsilon}
\newcommand{\eps}{\epsilon}

\newcommand{\Esymb}{{\bf E}}

\DeclareMathOperator*{\E}{\Esymb}

\usepackage[margin=2cm]{geometry}   

\newcommand{\RR}{\mathbb{R}}

\renewcommand{\[}{\begin{equation}}
\renewcommand{\]}{\end{equation}}

\definecolor{Blue}{rgb}{0,0,1}

\newcommand{\ba}{\mathbf{a}}

\providecommand{\RR}{\mathbb{R}}

\newcommand{\mc}[1]{\mathcal{#1}}

\newcommand{\bb}[1]{\mathbb{#1}}
\newcommand{\LF}{\mathrm{Lex}(\mc{F}, \mu)}
\newcommand{\tF}{\tilde{F}}

\begin{document}

\title{On the boundary of the region defined by homomorphism densities}


\author{Hamed Hatami \thanks{School of Computer Science, McGill University. \texttt{hatami@cs.mcgill.ca}. Supported by an NSERC grant.} \and Sergey Norin \thanks{Department of Mathematics and Statistics, McGill University. \texttt{snorin@math.mcgill.ca}. Supported by an NSERC grant.}}

\maketitle

\begin{abstract}
The Kruskal--Katona theorem together with a theorem of Razborov~\cite{MR2433944}  determine the closure of the set of points defined by the homomorphism density of the edge and the triangle  in finite graphs. The boundary of this region is a countable union of algebraic curves, and in particular, it is almost everywhere differentiable. One can more  generally consider the region defined by the homomorphism densities of a list of given graphs, and ask whether the boundary is as well-behaved as in the case of the triangle and the edge. Towards answering this question in the negative, we construct examples which show that the restrictions of the boundary to certain hyperplanes can have nowhere differentiable parts.
\end{abstract}

\section{Introduction}
For two simple finite graphs $F$ and $G$, \emph{the induced subgraph density} of $F$ in $G$, denoted by  $p(F;G)$, is  the number of induced copies of $F$ in $G$ divided by ${|V(G)| \choose  |V(F)|}$.  A large part of extremal graph theory is concerned with understanding the asymptotic relation between induced subgraph densities of different graphs.  In other words, given a list of graphs $\mc{F}=\{F_1,\ldots,F_k\}$, one would like to determine the compact region $T(\mc{F}) \subseteq [0,1]^k$ defined as the set of limit points of sequences $\{(p(F_1;G_i),\ldots,p(F_k;G_i))\}^{\infty}_{i=1}$, where $|V(G_i)| \to \infty$.  The general problem of understanding the  regions defined by subgraph densities was essentially first considered in~\cite{MR538044}. More recently, the study of graph limits~\cite{MR3012035} has brought new attention to such problems.

Describing the regions $T(\mc{F})$, even when $\mc{F}$ consists of a few small graphs, is a  major challenge. To the best of our knowledge the only non-trivial cases for which a complete description is known are the following. The set $T(K_2,K_3)$ has only been completely described less than a decade ago, and  Razborov's solution~\cite{MR2433944} to this problem is highly non-trivial, and is considered to be one of the major achievements of the theory of flag algebras~\cite{MR2371204}. More recently, Razborov's result was extended in two directions. Reiher~\cite{MR3549620} determined the inequalities defining the set $T(K_2,K_t)$ for every $t$, while Glebov et al.~\cite{GGHHKV} and Huang et. al~\cite{MR3200287} obtained the description of $T(\mc{F})$ for any pair $\mc{F}$ of graphs on three vertices.

The structure of $T(\mc{F})$ is very complex in full generality. In particular, it is shown in~\cite{MR2748400} that the separation problem for $T(\mc{F})$ is undecidable. That is there exists no algorithm which given $\mc{F}$ and an affine halfspace $S$ of $\bb{R}^k$ decides whether $T(\mc{F}) \cap S = \emptyset$ or not.

The focus of this paper is the boundary of the region $T(\mc{F})$. The boundary of $T(K_2,K_3)$ is a countable union of algebraic curves~\cite{MR2433944}, and, in particular,  is almost everywhere differentiable. This raises the question that whether the boundary is always as well-behaved; for example, is it almost everywhere differentiable?

Towards giving a strong negative answers to such questions we give an example which shows that the intersection of $T(\mc{F})$ with a hyperplane can have a nowhere differentiable boundary. The precise statement of our main result, Theorem~\ref{thm:main}, is technical and requires some preparation. In Section~\ref{sec:prel} we introduce the required terminology from the theory of graph limits. In Section~\ref{sec:lex} we prove some technical results, and finally in Section~\ref{sec:main} we state and prove our main result. We finish by some concluding remarks in Section~\ref{sec:conc}.

\section{Preliminaries \label{sec:prel}}

In this section we introduce the notions from the theory of graph limits that we will be using for the rest of the paper.

All graphs considered in this paper are simple and finite. For a graph $G$, let $V(G)$ and $E(G)$, respectively denote  the set of the vertices and the edges of $G$. The unique graph with no vertices is denoted by $\emptyset$.

\subsection{Weighted graphs and graphons}

We redefine $T(\mc{F})$ in terms of the induced homomorphism densities as follows.
For two graphs $F$ and $G$,  a map $\phi: V(F) \to V(G)$ is \emph{a strong homomorphism from $F$ to $G$} if $vw \in E(F)$ if and only if $\phi(v)\phi(w) \in E(G)$; that is $\phi$ preserves both adjacency and non-adjacency. Let $s(F,G)$ denote the number of  strong homomorphisms from $F$ to $G$, and let \emph{the induced homomorphism density}
$$t(F;G) := \frac{s(F,G)}{|V(G)|^{|V(F)|}},$$
be the probability that a random mapping (not necessarily injective) from the vertices of $F$ to the vertices of $G$ is a strong homomorphism.  We define $t(\emptyset;G) :=1$, for every graph $G$. For $\mc{F} = \{F_1,F_2,\ldots,F_k\}$, the set $T(\mc{F})$ defined in the introduction is (up to rescaling) equal to the closure of the set of points $(t(F_1;G), t(F_2;G), \ldots,t(F_k;G))$ for all $G$.

A \emph{weighted graph} $(G,\mu)$ is a graph $G$ together with a discrete probability measure $\mu$ on $V(G)$, such that $\mu(v)>0$ for every $v \in V(G)$. (The last assumption is unusual, but convenient.)
The quantity $t(F;G,\mu)$ can be defined similarly  to   $t(F;G)$ by mapping each vertex of $F$ to a vertex of  $G$  randomly and independently according to the distribution $\mu$, rather than uniformly.

A sequence of graphs $\{G_{i}\}_{i=1}^{\infty}$ is called \textit{convergent} if for every graph $F$, the sequence $\{t(F;G_{i})\}_{i=1}^{\infty}$ converges. Weighted graphs are not sufficient to represent the limits of all convergent sequences. However, Lov\'asz and Szegedy proved in~\cite{MR2274085} that the limit of a convergent graph sequence can be represented by a so-called \emph{graphon}, which is a symmetric measurable function $W\colon[0,1]^2 \to [0,1]$. More precisely, given a convergent graph sequence $\{G_n\}_{n=1}^\infty$, there always exists a graphon $W$ such that for every graph $F$, we have
\begin{equation}
 \label{eq:tW}
 \lim_{n \to \infty} t(F;G_{n}) =  \E \left[\prod_{\substack{(u,v)\in E(F)}} W(x_{u},x_{v})\prod_{(u,v) \notin E(F)} (1 - W(x_{u},x_{v}))\right],
\end{equation}
where $\{x_{u} \ | \ u \in V(F)\}$ are independent random variables  taking values uniformly in $[0,1]$.
We denote the expected value in the right-hand side by $t(F;W)$. Conversely, for every graphon $W$, one can construct a  graph sequence that converges to $W$ in the above sense.
Note that a convergent sequence of graphs does not converge to a \emph{unique} graphon. We say that graphons $W$ and $W'$ are \emph{weakly isomorphic} if $t(F;W)=t(F;W')$ for every graph $F$. Thus a convergent sequence of graphs does converge to a unique graphon up to weak isomorphism. It will be convenient for us not to distinguish between weakly isomorphic graphons in this paper.

For every weighted graph $(G,\mu)$ (on vertex set $[n]$), we define a step-function $W_{G,\mu}$ as follows: split $[0,1]$ into $n$ intervals $J_1,\ldots,J_n$ of lengths $\lambda(J_i)=\mu(i)$, and for $x \in J_i$ and $y \in J_j$, set $W_{G,\mu}(x,y)$ to $1$, if $i$ is adjacent to $j$ in $G$, and to $0$, otherwise. Note that we have $t(F;G,\mu) = t(F;W_{G,\mu})$ for every  graph $F$ and weighted graph $(G,\mu)$. In this sense graphons are a natural extension of weighted graphs (and, in particular, graphs).

\subsection{Quantum graphs}\label{sec:quantum}
Let $L$ be a finite (possibly empty) subset of natural numbers. An \emph{$L$-labeled graph} is a graph in which some of the vertices are labeled by distinct elements of $L$, such that every label is assigned to exactly one vertex (there may be any number of unlabeled vertices). Let ${\mathcal H}_L$ denote the set of all $L$-labeled graphs up to label-preserving isomorphism. For brevity we denote the set of unlabeled graphs $\mc{H}_{\emptyset}$ by $\mc{H}$. A graph is \emph{partially labeled} if it is $L$-labeled for some set $L$.  A partially labeled graph in which all vertices are labeled is called a \emph{fully labeled} graph.

We extend the definition of the induced homomorphism density to partially labeled graphs in the following way. Consider a finite set $L \subset \mathbb{N}$, a partially labeled graph $H \in {\mathcal H}_L$, a graph $G$, and a map $\phi\colon L \to V(G)$. Then $t(H, \phi;G)$ is defined to be the probability that a random map from $V(H)$ to $V(G)$ is a strong homomorphism conditioned on the event that the labeled vertices are mapped according to $\phi$.

The definition of $t(H, \phi;G)$ extends to graphons in a straightforward manner. That is, given a graphon $W$ and a map $\phi\colon L \to [0,1]$, the induced homomorphism density $t(H,\phi;W)$ is the expected value  in (\ref{eq:tW}) but now for every $i \in L$, we fix $x_{u_i} = \phi(i)$ where $u_i$ is the vertex with label $i$, and the randomness is over the rest of the variables.

Consider real numbers $\alpha_1,\ldots,\alpha_k$ and graphs $H_1, \ldots, H_k$. In order to study the linear combinations of induced homomorphism densities, it is useful to define quantum graphs and labeled quantum graphs.  An \emph{$L$-labeled quantum graph} is an element of the vector space $\mathbb{R}[{\mathcal H}_L]$, i.e. it is a formal linear combination of graphs in $\mc{F}_L$. For a graph $G$, an $L$-labeled quantum graph, $f = \sum_{i=1}^n \alpha_i H_i \in \mathbb{R}[{\mathcal H}_L]$  and a map $\phi\colon \cup_{i=1}^k L_i \to V(G)$, define $t(f, \phi;G) := \sum_{i=1}^k \alpha_i t(H_i,\phi|_{L_i} ;G)$.

We  define the product $H_1 \cdot H_2$ of two vertex disjoint $L$-labeled graphs $H_1$ and $H_2$, to be zero if the labeled subgraphs of $H_1$ and $H_2$ induced by the labeled sets differ, and, otherwise, we define $H_1 \cdot H_2$ to be the sum of all distinct $L$-labeled graphs obtained by identifying the vertices of $H_1$ and $H_2$ with the same labels and possibly adding some edges between unlabeled vertices of $H_1$ and unlabeled vertices of $H_2$.  This product is defined so that for every two partially labeled graphs $H_1,H_2 \in \mc{F}_L$, a graph $G$, and a map $\phi\colon L \to V(G)$,
\begin{equation}
\label{eq:product}  t(H_1 \cdot H_2, \phi; G) = t(H_1, \phi;G) t(H_2, \phi; G).
\end{equation}
We extend this product to $\mathbb{R}[\mc{F}_L]$ by linearity; that is if $f = \sum_{i=1}^n \alpha_i H_i \in \mathbb{R}[\mc{F}_L]$ and $g = \sum_{i=1}^m \beta_i F_i \in \mathbb{R}[\mc{F}_L]$, then $f \cdot g =  \sum_{i=1}^n \sum_{j=1}^m \alpha_i \beta_j  H_i \cdot F_j$. Consider a finite set $L \subset \mathbb{N}$, a graph $G$, and a map $\phi\colon L \rightarrow V(G)$.
It follows from (\ref{eq:product}) that $f \mapsto t(f, \phi;G)$ defines a homomorphism from  $\mathbb{R}[\mc{F}_L]$ to $\RR$.

Let the linear map $\llbracket  \cdot \rrbracket: \bigcup_L \mathbb{R}[\mc{F}_L] \rightarrow \mathbb{R}[{\mathcal F}]$ be defined by un-labeling all the labeled vertices. Note that $t(\llbracket  f^2 \rrbracket,W) = 0$ for a graphon $W$ and an $L$-labeled quantim graph $f \in \mathbb{R}[\mc{F}_L]$ if and only if  $t(f,\phi;W)=0$ for almost every $\phi:L \to [0,1]$.

\subsection{Infinite lexicographic products}

The lexicographic product $F \otimes  H$ of graphs $F$ and $H$ is the graph with vertex set $V(F) \times V(H)$ where two vertices $(u,v)$ and $(x,y)$ are adjacent in $F \otimes H$ if and only if either $u$ is adjacent to  $x$ in $F$, or $u = x$ and $v$ is adjacent to $y$ in $H$. The lexicographic product of two weighted graphs is defined by $$(F,\mu) \otimes (H,\nu) = (F \otimes H, \mu \times \nu).$$
Given an infinite sequence of weighted graphs $(F_i,\mu_i)_{i \in \mathbb{N}}$, \emph{the (infinite) lexicographic product $\otimes_{i=1}^\infty (F_i,\mu_i)$} is the graphon defined (up to weak isomorphism) as the limit of the convergent sequence of weighted graphs $\{\otimes_{i=1}^n (F_i,\mu_i)\}_{n \in \bb{N}}$. If all the graphs in the sequence  $(F_i,\mu_i)_{i \in \mathbb{N}}$ are equal to the weighted graph $(F,\mu)$ then the resulting product is called \emph{the infinite lexicographic power of $(F,\mu)$} and is denoted by $\otimes^{\infty}(F,\mu)$.

Let $W$ be a graphon, and let $(F,\mu)$ be a weighted graph. We say that a partition $(J_v)_{v \in V(F)}$ of the interval $[0,1]$ is an \emph{$(F,\mu)$-partition of $W$} if
\begin{itemize}
	\item $\lambda(J_v) = \mu(v)$ for all $v \in V(F)$, and
	\item for every pair of distinct vertices $u,v \in V(F)$, the restriction of $W$ to $J_u \times J_v$ is identically $1$, if $u$ and $v$ are adjacent in $F$, and is identically $0$ otherwise.
\end{itemize} Note that this implies no restrictions on the value of $W$ on the diagonal cells $J_u \times J_u$ for $u \in V(G)$.
We say that a graphon  $W$ is \emph{$(F,\mu)$-partitionable} if it admits an  $(F,\mu)$-partition.

For a graphon $W$ and a measurable set $J\subseteq [0,1]$ we define the \emph{restriction $W[J]:[0,1]^2 \to [0,1]$ of $W$ to $J$} by taking any measurable bijection $\phi: [0,1] \to J$ such that $\lambda(S)=\lambda(\phi^{-1}(S))\lambda(J)$ for every measurable $S \subseteq J$, where $\lambda$ denotes the Lebesgue measure, and defining  $W[J](x,y)=W(\phi(x),\phi(y))$. The choice of $\phi$ is irrelevant for our purposes as the graphons obtained using different choices are weakly isomorphic to each other.

Given an  $(F,\mu)$-partition $(J_v)_{v \in V(F)}$ of a graphon $W$, let $W_v=W[J_v]$ for $v \in V(G)$.
We say that $(W_v)_{v \in V(F)}$ is the collection of the \emph{bags} of the  $(F,\mu)$-partition $(J_v)_{v \in V(F)}$.

It is important for the application of the above definitions to note that for any infinite lexicographic product  $W=\otimes_{i=1}^\infty (F_i,\mu_i)$ of weighted graphs and every $n \geq 1$, there exists a  $\otimes_{i=1}^n (F_i,\mu_i)$-partition of $W$ such that all the bags of the partition are weakly isomorphic to   $\otimes_{i=n+1}^\infty (F_i,\mu_i)$.

\subsection{Blowups}\label{s:blowup}
In addition to lexicographic products, in some of the technical arguments we will use another, simpler, construction, called blowups. For a graph $F$ and  a positive integer vector $\ba \in \mathbb{N}^{V(F)}$, the $\ba$-\emph{blow-up} of $F$, denoted by $F^{(\ba)}$, is the graph obtained by replacing every vertex $v$ of $H$ with $\ba(v)$ different vertices where a copy of $u$ is adjacent to a copy of $v$ in the blow-up graph if and only if $u$ is adjacent to $v$ in $F$. If $\ba$ is equal to $d$ in every coordinate we will  denote $\ba$-blow-up of $F$ by $F^{(d)}$.
	
If $F$ is a fully labeled graph then the blow up  $F^{(\ba)}$ is defined in a similar manner, where the newly added vertices will be unlabeled and the original vertices of $F$ will keep their labels. Let $\tilde{F}^{(\ba)}$ be the quantum graph equal to the sum over all graphs that can be obtained from $F^{(\ba)}$ by possibly adding some edges between different copies of the same vertex of $F$.

Finally, for $I \subseteq V(F)$, let $F^{I}$ denote $F^{(\ba_I)}$ where  $\ba_I(v)=2$ if $v \in I$, and $\ba_I(v)=1$, otherwise. That is, $F^{I}$ is obtained from $F$ by replacing every vertex in $I$ by a pair of twin vertices. Similarly $\tilde{F}^{I}$ is defined  as  $\tilde{F}^{(\ba_I)}$.

\subsection{Finite forcibility}\label{s:forcibility}

A  graphon $W$ is \emph{finitely forcible}  if there exists a quantum graph $f$ such that $t(f;W')=0$ for a graphon $W'$ if and only if $W'$ is weakly isomorphic to $W$.  Finitely forcible graphons, first formally defined in~\cite{MR2802882},  are  the structures that can appear as ``unique" solutions of extremal problems involving graph homorphism densities. As mentioned in the introduction, such problems can be extremely difficult in full generality. Increasingly complex finitely forcible graphons are constructed in~\cite{MR3279390,GKKFin,MR2802882}, furthering the evidence in support of this claim. In particular, it is asked in~\cite{MR2802882} for which graphs $F$  the infinite lexicographic power $\otimes^{\infty}F$ is finitely forcible. Our main technical result, Lemma~\ref{lem:key}, provides a partial answer to this question.

Extending the above definition, we say that a family $\mc{W}$ of graphons is \emph{finitely forcible} if there exists a quantum graph $f$ such that $t(f;W')=0$ for a graphon $W'$ if and only if $W'$ is weakly isomorphic to some graphon $W \in \mc{W}$. We say that such a quantum graph $f$ \emph{forces} $\mc{W}$.

\section{Forcing lexicographic products}\label{sec:lex}

In this section we prove a series of technical results related to finite forcibility, which are used in the proof of our main result in the next section.

Let $H$ and $G$ be graphs. We say that a map $\phi:V(G) \to V(H)$ is \emph{a folding of $G$ into $H$} if for every pair of vertices $u,v \in V(G)$, either $\phi(u)=\phi(v)$, or $uv \in E(G)$ if and only if $\phi(u)\phi(v) \in E(H)$. That is $\phi$ preserves adjacencies and non-adjacencies except that it does not restrict adjacencies of the vertices of $G$ mapped to the same vertex of $H$. Thus every strong homomorphism is a folding, but not vice versa. We say that a folding $\phi$ is \emph{trivial} if it maps all the vertices to a single vertex.

A set $A \subseteq V(H)$ is called \emph{homogeneous} in a graph $H$, if for every pair of distinct vertices $u,v \in A$, $N(u) \backslash A = N(v) \backslash A$, where $N(u)$ and $N(v)$ respectively denote the set of the neighbors of $u$ and $v$. Equivalently, every vertex outside $A$ is either adjacent to all the vertices in $A$, or has no neighbor in $A$.  We call a graph $H$ \emph{prime}, if it does not contain any homogeneous sets $A$ with $1<|A| \le |V(H)|-1$. Note the following easy property of prime graphs.

\begin{proposition}\label{p:foldinggraph}	Let  $\phi: V(H) \to V(G)$ be a folding of a graph $H$ into a graph $G$. Let $A$ be a homogeneous set in $G$. Then $\phi^{-1}(A)$ is homogeneous in $H$. 	
In particular, if $H$ is prime, then either $\phi$ is trivial, or $\phi$ is injective, and thus is a strong homomorphism.
\end{proposition}	

We say that a graph $H$ is \emph{stringent} if $H$ is prime, and  furthermore, $H$ does not have any non-identity automorphisms. If a graph $H$ is stringent, then in particular the identity map is the only map from $H$ to itself that preserves both adjacency and non-adjacency.

\begin{proposition}
\label{prop:randomStringent}
For any fixed $p \in (0,1)$ the Erd\"os-Renyi  random graph $G(n,p)$ is asymptotically almost surely stringent.
\end{proposition}
\begin{proof}
Erd\H{o}s and Renyi~\cite{MR0156334} have shown that the automosphism group of $G(n,p)$ is a.a.s. trivial. Thus it remains to show that $G(n,p)$ is a.a.s. prime.

If $G=G(n,p)$ is not prime  then there exists a set $S$ of size $2 \le k \le n-1$ such that every vertex in $V(G)- S$ is either adjacent to all vertices in $S$ or is not adjacent to any vertex in $S$. Thus the probability that $G$ is not prime is at most
$$\sum_{k=2}^n {n \choose k} \left((1-p)^k+p^k\right)^{n-k} \leq \sum_{k=0}^n {n \choose k} \left((1-p)^2+p^2\right)^{n-k} = \left((1-p)^2+p^2\right)^n \underset{n \to \infty}{\longrightarrow} \:0,$$
as desired.
\end{proof}

We will need a statement similar to Proposition~\ref{p:foldinggraph} involving maps from a prime graph into a graphon. Extending the definition of homogeneous sets in graphs,
we say that a set $S \subseteq [0,1]$ is \emph{weakly homogeneous} in a graphon $W$, if there exists a function $e:[0,1] \setminus S \to \{0,1\}$ such that $W(x,y)=e(x)$ for almost every pair $(x,y)$ with $x \in [0,1] \setminus S$ and $y \in S$. We say that  $S \subseteq [0,1]$ is \emph{strongly  homogeneous} instead if  $W(x,y)=e(x)$ holds for all pairs $(x,y)$ as above.
As we do not distinguish between weakly isomorphic graphons, and, in particular, between graphons which differ on a set of measure zero, it will be  convenient for us to assume that a weakly homogeneous set is strongly homogeneous.

We are now ready to state a graphon analogue of Proposition~\ref{p:foldinggraph}.

\begin{proposition}\label{p:foldinggraphon}
Let $W$ be a graphon, let $S$ be a strongly homogeneous set in $W$,  and let $H$ be a fully labelled graph. If $t(H,\phi;W)> 0$ for an injective  $\phi:V(H) \to [0,1]$, then  $\phi^{-1}(S)$ is homogeneous in $H$. In particular, if $H$ is prime, then
 either $\phi^{-1}(S)=V(H)$, or $|\phi^{-1}(S)| \leq 1$.
\end{proposition}

Note that if $(J_v)_{v \in V(F)}$ is an  $(F,\mu)$-partition  of a graphon $W$, then, in particular, the set $J_v$ is strongly homogeneous for every $v \in V(F)$. For a measure $\mu$ on the vertex set of a graph $G$ and a positive integer $k$, define $\mathbf{m}_k(\mu)=\sum_{v \in V(G)}\mu^k(v)$. The preceding observation together with Proposition~\ref{p:foldinggraphon} implies the following.

\begin{lemma}\label{lem:density}
	Let $(F,\mu)$ be a weighted graph, let $W$ be an $(F,\mu)$-partitionable graphon, and let $(W_v)_{v \in V(F)}$ be the collection of bags of some $(F,\mu)$-partition of $W$. Then
	\begin{equation}\label{e:density2}t(H;W)=t(H;F,\mu)+\sum_{v \in V(F)}\mu^k(v)t(H;W_v)
	\end{equation}
	for any prime graph $H$ on $k$ vertices.
		
	In particular, if $W=\otimes_{i=1}^\infty (F_i,\mu_i)$ for some sequence of weighted graphs $ (F_i,\mu_i)_{i \in \bb{N}}$, then
		\begin{equation}\label{e:density}
		t(H;W) = \sum_{i=1}^\infty \left( \prod_{j=1}^{i-1}\mathbf{m}_k(\mu_j)\right) t(H;F_i,\mu_i).
		\end{equation}
\end{lemma}	

 We say that a finite collection $\mc{F}$  of  graphs is \emph{$n$-stringent} if each graph $F \in \mc{F}$ is stringent, and $V(F)=[n]$ for every  $F \in \mc{F}$. Let $\mu: [n] \to (0,1)$ be a probabilistic measure.  We define $\mathrm{Lex}(\mc{F},\mu)$ to  be the set of all infinite lexicographic products  $\otimes_{i=1}^\infty J_i$, where for every $i \ge 1$ we have $J_i = (F, \mu)$  for some $F \in \mc{F}$. The following lemma is the main result of this section.

\begin{lemma}
	\label{lem:key}
	The family $\mathrm{Lex}(\mc{F}, \mu)$ is finitely forcible for every $n$-stringent collection of graphs $\mc{F}$ and every probabilistic measure $\mu: [n] \to (0,1)$.
\end{lemma}

The remainder of the section is occupied with the proof of Lemma~\ref{lem:key}. The proof relies, in particular, on the following simple observation about finitely forcible families.

\begin{proposition}\label{p:intersection}
	Let $\mc{W}_1,\mc{W}_2$ be finitely forcible families of graphons closed under weak isomorphism. Then $\mc{W}_1 \cap \mc{W}_2$ is finitely forcible.
\end{proposition}
\begin{proof}
	If $f_i$ forces $\mc{W}_i$ for $i=1,2$, then $f_1^2+f_2^2$ forces $\mc{W}_1 \cap \mc{W}_2$.
\end{proof}

A graphon is called \emph{random-free} if it is  $\{0, 1\}$-valued almost everywhere.

\begin{proposition}
\label{prop:randomfree}
Let $\mc{F}$ be a family of $n$-stringent graphs. Then $\mathrm{Lex}(\mc{F}, \mu) \subseteq \mc{W}_{\mathrm{rf}}$ for a finitely forcible family of random-free graphons $\mc{W}_{\mathrm{rf}}$.
\end{proposition}
\begin{proof}
 Consider a bipartite graph $H$ with a bipartition $V(H)=X \cup Y$, and let the quantum graph $\overline{H}$ be the sum of all graphs that can be obtained from $H$ by adding only edges whose end-points are either both in $X$ or both in $Y$. Note that $t(\overline{H};G)$ is the probability that a random map $\phi:V(H) \to V(G)$ satisfies $\phi(x)\phi(y)\ \in E(G)$ if an $xy \in E(H)$ for all pairs $x \in X$, $y \in Y$. It is shown in~\cite{MR2815610}  that $t(\overline{H};W)=0$ implies that $W$ is random-free. Thus our task reduces to finding a bipartite graph $H$ such that $t(\overline{H};W)=0$ for all $W \in \mathrm{Lex}(\mc{F}, \mu)$.

 It is well-known that with high probability the size of the largest independent set in the Erd\"os-Renyi  random graph $G(m,1/2)$ is bounded by $2 \log m$. Thus by Proposition~\ref{prop:randomStringent}, there exists a stringent graph $K$ with $m:=|V(K)|>n$ whose independence number is smaller than $2 \log m$. Let $H$ be the bipartite graph obtained from $K$ by splitting every vertex $i$ of $K$ into two vertices $u_i$ and $v_i$, and connecting  $u_i$ to $v_j$ if and only if $ij \in E(K)$. The sets $X=\{u_i:i \in V(K)\}$ and  $Y=\{v_i:i \in V(K)\}$ form a bipartition of $H$.  Suppose $t(\overline{H};W)>0$ for some  $W \in \mathrm{Lex}(\mc{F}, \mu)$. Then for some $F \in \mathcal{F}$, there exists a map $\phi:V(H) \to V(F)$ such that $t(\overline{H},\phi;F,\mu)>0$, and $\phi$ does not map all the vertices of $H$ to a single vertex. It follows from  $t(\overline{H},\phi;F,\mu)>0$ that for all $u_i \in X$ and $v_j \in Y$ satisfying $\phi(u_i) \neq \phi(v_j)$, we have $\phi(u_i)\phi(v_j)\ \in E(F)$ if and only if $u_iv_j \in E(H)$.
 Pick $a \in V(F)$ such that $2m/n \le |\phi^{-1}(a)|<2m$.  Note that $\{i:  \phi(u_i) = \phi(v_i) = a\}$ is  a homogeneous set  in $K$, and thus is of size at most $1$. But both sets $\{i:  \phi(u_i) = a, \phi(v_i) \neq a\}$ and $\{i:  \phi(u_i) \neq a, \phi(v_i) = a\}$ are independent sets in $K$, and thus they are of size at most $2 \log m$. For sufficiently large $m$, this contradicts the choice of $a$.
 \end{proof}

By Proposition~\ref{p:intersection}, to prove Lemma~\ref{lem:key} it suffices to construct finitely forcible families of graphons $\mc{W}_1,\ldots,\mc{W}_l$ for some $l$ such that $\LF = \cap_{i=1}^l \mc{W}_i$. Defining these families will require a technical definition.

\begin{definition}
\label{def:splitGraphons}
Let $k$ be a positive integer, $(F,\mu)$ a weighted graph, and $W$ a graphon. We say that $W$ is an \emph{$(F,\mu,k)$-split graphon}  if for (almost) every  map $\phi:V(F) \to (0,1)$ with $t(F,\phi;W)>0$, there exists a homogeneous set $S$, and a partition $(J_v)_{v \in V(F)}$ of $S$ such that the following holds:
\begin{description}
	\item[(S1)] $(J_v)_{v \in V(F)}$  corresponds to an $(F,\mu)$-partition of $W[S]$;
	\item[(S2)] $\phi(v) \in J_v$ for every $v \in V(F)$;
	\item[(S3)] $t(H;W[J_v])=t(H;W[J_u])$ for every pair of vertices $u,v\in V(F)$, and every graph  $H$ on $k$ vertices.
\end{description}
\end{definition}

Let $\mc{S}_{k}(F,\mu)$ denote the family of all $(F,\mu,k)$-split graphons.

\begin{lemma}\label{l:split}
	Let $k$ be a positive integer, and let $F$ be a stringent graph. Let  $\mc{W}_{\mathrm{rf}}$ be a finitely forcible family of random-free graphons. Then  $\mc{S}_{k}(F,\mu) \cap  \mc{W}_{\mathrm{rf}}$ is finitely forcible for every probabilistic measure $\mu:V(F) \to (0,1)$.
\end{lemma}

\begin{proof}
By Proposition~\ref{prop:randomfree} it suffices to show the existence of a quantum graph $f$ such that a random-free graphon $W$ satisfies $t(f;W)=0$ if and only if $W \in \mc{S}_{k}(F,\mu)$.
Assume  $V(F)=[n]$. We will think of $F$ as an $[n]$-labelled graph.

Consider a map $\phi:[n] \to [0,1]$ and a $\{0,1\}$-valued graphon $W$ such that $t(F,\phi;W)=1$. For every vertex $i$ of $F$, let $J_i$ be the set of points $x \in [0,1]$ such that $W(x,\phi(j))=W(\phi(i),\phi(j))$  for all $j \in [n]$.  Let $J_\circ$  be the set of $x \in [0,1]$ such that $W(x,\phi(j))=0$ for all $j \in [n]$, and let $J_\bullet$ be the set of $x \in [0,1]$ such that $W(x,\phi(j))=1$ for all $j \in [n]$. Since $F$ is stringent, the sets $J_1,\ldots,J_n,J_\circ,J_\bullet$ are disjoint.  Let $S:=J_1 \cup \ldots \cup J_n$. We would like to force the sets $S$ and $J_1,\ldots,J_n$ to satisfy the conditions of Definition~\ref{def:splitGraphons}. We will introduce a few conditions of the form $t(f_i,\phi;W)=a_i$ that would imply this. Consequently by the discussion in the last paragraph of Section~\ref{sec:quantum}, setting $f=\llbracket \sum (f_i-a_i)^2 \rrbracket$ will yield the desired result.

Define $h =\sum_{1 \leq i}^{n} \tF^{\{i\}}$,  let $F^{\circ}$  be obtained from $F$ by adding an unlabeled isolated vertex, and let $F^{\bullet}$ be obtained from $F$ by connecting a new unlabeled vertex to all the vertices of $F$. For every graph $H$ on  $k$ vertices, and every $i \in [n]$, let $F \oplus_i H$ be the $[n]$-partially labelled graph that is obtained from $F$ by first creating $k$ unlabeled twins  for the vertex $i$, and then implanting a copy of $H$ on these unlabeled vertices. Let $F \tilde{\oplus_i} H$ be sum of all the $2^{k}$ graphs that are obtained from $F \oplus_i H$ by possibly adding some edges between the vertex $i$ and the unlabeled vertices. Note the following
\begin{enumerate}
 \item The condition  $t(h+F^\circ+F^\bullet,\phi;W)=1$ implies that $J_1 \cup \ldots \cup J_n \cup J_\circ \cup J_\bullet=[0,1]$ up to a measure zero difference.
 \item The condition $t(\tF^{\{i\}} - \mu(i)h,\phi;W)=0$ implies that $\lambda(J_i)=\mu(i) \lambda(J_1 \cup \ldots \cup J_n)$, and thus $J_i$ are of desired measures.
 \item The condition $t(\tF^{\{i,j\}} - \tF^{\{i\}}\tF^{\{j\}},\phi;W)=0$ implies that for almost every $x \in J_i$ and $y \in J_j$, we have $W(x,y)=1$ if and only if $ij \in E(F)$. This together with the previous condition imply that  $(J_i)_{i \in V(F)}$ is an $(F,\mu)$-partition of $W[S]$.
 \item It can be similarly forced that $W(x,y)=1$ for all $x \in J_i$ and $y \in J_\bullet$, and $W(x,y)=0$ for all $x \in J_i$ and $y \in J_\circ$. Thus $S$ is a homogeneous set.
 \item The condition $\mu(i)^{-k}t( F \tilde{\oplus_i} H,\phi;W)=\mu(j)^{-k} t( F \tilde{\oplus_j} H,\phi;W)$ implies  $t(H;W[J_i])=t(H;W[J_j])$, and thus (S3) can be forced in this manner.
\end{enumerate}
As we mentioned above, the desired quantum graph can be constructed from these conditions.
\end{proof}

Before finally stating the proof of Lemma~\ref{lem:key}, let us make a simple observations regarding $(F,\mu,k)$-split graphons.

\begin{lemma}
\label{l:observ}
Let $F_1$ and $F_2$ be stringent graphs, and let $W$ be an $(F_1,\mu,k)$-split graphon. Let $S$, and $(J_v)_{v \in V(F_1)}$ be as in Definition~\ref{def:splitGraphons} for some $\phi:V(F_1) \to [0,1]$. Then the following hold.
\begin{description}
 \item [(i)] If $W$ is an $(F_2,\mu,k)$-split graphon, then so is $W[J_v]$.
 \item [(ii)] If $\phi:V(F_2^{(d)}) \to [0,1]$ satisfies $t(\tF_2^{(d)},\phi;W)>0$, then either $\phi(V(F_2^{(d)})) \subseteq S$, or there is a vertex $v \in V(F_2)$ such that for all the vertices $u \in V(F_2^{(d)}) $ that are not copies of $v$, we have $\phi(u) \not\in S$.
\end{description}
\end{lemma}
\begin{proof}
To verify (i) pick a vertex $v_0 \in V(F_1)$, and consider a map $\psi$ with $t(F_2,\psi;W[J_{v_0}])>0$ and its corresponding map $\psi':V(F_i) \to J_{v_0}$ with $t(F_2,\psi';W)>0$.  If $W$ is an $(F_2,\mu,k)$-graphon, one can find a homogeneous $\tilde{S}$ in $W$ and a partition $(\tilde{J}_v)_{v \in V(F)}$ of $\tilde{S}$ satisfying (S1), (S2), and (S3). Since $F_2$ is stringent and $J_{v_0}$ is a homogeneous set for $W$, by Proposition~\ref{p:foldinggraphon} we must have $\tilde{J}_v \subseteq J_{v_0}$ for all $v$.  Thus $W[J_{v_0}]$ is an $(F_2,\mu,k)$-split graphon.

To verify (ii), consider a map $\phi:V(F_2^{(d)}) \to [0,1]$ satisfying $t(\tilde{F_2}^{(d)},\phi;W)>0$. Consider a copy of $F_2$ in $F_2^{(d)}$. If $\phi$ maps two different vertices from this copy to points in $S$, then since $S$ is homogeneous and $F_2$ is stringent, it must map all the vertices of this copy into $S$. This in turn implies that all the vertices of  $F_2^{(d)}$ must be mapped into $S$.
\end{proof}

\begin{proof}[Proof of Lemma~\ref{lem:key}]	
  Let $\mc{F}=(F_i)_{i=1}^l$ be the family of $n$-stringent graphs; that is for every $i$, the graph $F_i$ is stringent and $V(F_i)=[n]$. We consider the graphs $F_1,\ldots,F_n$ as fully labeled graphs. Let $d$ be a positive integer chosen to be sufficiently large to satisfy the inequalities appearing later in the proof. By Lemma~\ref{l:split}, the class of graphons $\mc{W}^*=\cap_{i=1}^l \mc{S}_{nd}(F_i,\mu) \cap \mc{W}_{\mathrm{rf}}$ is finitely forcible, where $\mc{W}_{\mathrm{rf}}$ is as in Proposition~\ref{prop:randomfree}. Clearly $\LF \subseteq \mc{W}^*$. Let $f_1$ be a quantum graph that forces $\mc{W}^*$. We will show that there exists a quantum graph $f_2$ such that $t(f_2;W)=0$ for $W \in \mc{W}^*$ if and only if $W \in \LF$. This will imply that $f_1^2+f_2^2$ forces $\LF$, and thus the lemma.

	Let $f_2'=\llbracket \sum_{i=1}^l\tilde{F}_i^{(d)} \rrbracket$, and set $\rho =\prod_{i=1}^n\mu(i)$ and $m=\frac{\rho^d}{1-\mathbf{m}_{nd}(\mu)}$. By (\ref{e:density}) we have $t(f_2';W)=m$ for every $W \in \LF$. We will show that, if $t(f_2';W) \geq m$ for some $W \in \mc{W}^*$, then  $W \in \mc{W}^*$. This will imply that the quantum graph $f_2=f'_2-m$ is as desired.
	
	Let $m^*=\max_{W \in  \mc{W}^*}t(f_2';W)$. Our first goal is to show that $m^*=m$. Consider $W \in  \mc{W}^*$  with $t(f_2';W) = m^*$, and pick $i \in [l]$ so that  $t(\llbracket \tilde{F}_i^{(d)} \rrbracket;W) \geq m/l$. Consequently, there exists $\phi: [n] \to [0,1]$ with 	$t(\tilde{F}_{i}^{(d)},\phi;W) \ge m/l$;  in particular, $\phi$ satisfies the conditions in the definition of the $(F_i,\mu,nd)$-split graphon. Consider the homogeneous set $S$, and the partition $(J_v)_{v \in [n]}$ as in Definition~\ref{def:splitGraphons}. Note that $W[S]  \in \mc{W}^*$, by construction, and moreover since $F_1,\ldots,F_l$ are stringent graphs, by  Lemma~\ref{l:observ}~(i)  for every $v \in [n]$, $W[J_v] \in \mc{W}^*$.

	Set $\epsilon=1-\lambda(S)$, and note that
	$$
 	\frac{m}{l} \leq t(\tilde{F}_{i}^{(d)},\phi;W)  \leq (1-\eps)^{n(d-1)}.
	$$
	Therefore,
	\begin{equation}\label{e:bound1}
	\eps \leq 1-\left(\frac{m}{l}\right)^{\frac{1}{n(d-1)}} = 1- \left(\frac{\rho^d}{l(1-\mathbf{m}_{nd}(\mu))}\right)^{\frac{1}{n(d-1)}} \leq \frac{1}{2},
	\end{equation}
	provided that $d$ is sufficiently large. However, by  Lemma~\ref{l:observ}~(ii), we have
	\begin{align}\label{e:bound2}
	m^* &= t(f_2';W) \leq (1-\eps)^{nd}t(f_2'; W[S])+l(n+1) \eps^{n(d-1)} \notag
	\\ &\leq (1-\eps)^{nd}m^*+l(n+1) \eps^{n(d-1)} \leq (1-\eps)m^* + l(n+1) \eps^{d-1}.
	\end{align}
	If $\eps>0$, then (\ref{e:bound1}) and (\ref{e:bound2}) imply that
	$$m \leq m^* \leq l(n+1)\eps^{d-2} \leq \frac{l(n+1)}{2^{d-2}};$$
	a contradiction for sufficiently large $d$.  Hence we conclude that $\eps=0$. Thus $W$ is $F_{i}$-partitionable, and by Lemma~\ref{lem:density},
    \begin{align}\label{e:bound3}
    m^* &= t(f_2'; W)= \rho^d + \sum_{v \in [n]}\mu(v)^{nd} t(f_2';W[J_v]) \leq  \rho^d + \sum_{v \in [n]}\mu(v)^{nd} m^*=\rho^d + \mathbf{m}_{nd}(\mu) m^*.
    \end{align}
    It follows from the definition of $m$ that $m^*=m$, and that all the inequalities in  (\ref{e:bound3}) hold with equality; that is $t(f_2';W[J_v])=m$ for all $v \in [n]$. Thus, applying the preceding argument to $W[J_v]$ for every $v\in [n]$, we conclude that each $W[J_v]$ is $(F_{j(v)},\mu)$-partitionable for some $1 \leq j(v) \leq l$.

    Note that $t(\llbracket \tilde{F}^{(d)}_{j(v)} \rrbracket;W[J_v]) \geq \rho^d$. Moreover, by  Lemma~\ref{lem:density}, we obtain that if $k \neq j(v)$, then
    $$t(\llbracket \tilde{F}^{(d)}_{k}\rrbracket;W[J_v]) \leq \sum_{u \in [n]}\mu(u)^{nd} m = \mathbf{m}_{nd}(\mu) m < \rho^d,$$
    as long as $d$ is large enough so that $\mathbf{m}_{nd}(\mu)< \frac{1}{2}.$
    Thus by the property (S3) of $(F_i,\mu,nd)$-split graphons we deduce that there exists $j$ such that $W[J_v]$ is $F_j$-partitionable for every $v \in [n]$. Repeating the same argument recursively we deduce $W \in \LF$, as desired.
\end{proof}

\section{The main result}\label{sec:main}

We are now ready to state and prove our main result, which was informally described in the introduction, using the tools developed in Section~\ref{sec:lex}.

\begin{theorem}
\label{thm:main}
There exists a quantum graph $f$ such that the function $\psi \colon[0,1] \to \mathbb{R}$ defined as $\psi(x)=\min t(K_2;W)$  where the minimum is over all graphons $W$ with $t(f;W)=0$ and $t(P_4; W) = x$ is nowhere differentiable in the interval $[\alpha_1,\alpha_2]$ for some $0 \le \alpha_1 < \alpha_2 \le 1$.
\end{theorem}

\begin{proof}
Let $F_1$ and $F_2$ be two stringent graphs with $V(F_1)=V(F_2)=[n]$. Let $\mu$ be the unique probability measure on $[n]$ satisfying $\sum_{i=1}^n \mu(i)^4 =\frac{1}{2}$ and $\mu(2)=\ldots=\mu(n)$. Let $\alpha_i = t(P_4; F_i,\mu)$, and $\beta_i= \frac{t(P_2; F_i,\mu)}{1-\mathbf{m}_2(\mu)}$.  We assume   $\alpha_1 < \alpha_2$, and $\beta_1 < \beta_2$. Indeed if one takes $F_1$ to be a random graph $G(n,1/3)$, and  $F_2$ to be a random graph sampled $G(n,1/2)$, then it is straightforward to see, using Proposition, \ref{prop:randomStringent} that $F_1$ and $F_2$  satisfy the above properties asymptotically almost surely.

Let $\mc{F} =\{(F_1,\mu),(F_2,\mu)\}$. Then $\mc{F}$ is stringent, and thus by Lemma~\ref{lem:key} there exists a quantum graph $f$ such that $t(f;W')=0$ for a graphon $W'$ if and only if $W'$ is weakly isomorphic to some graphon $W \in \mathrm{Lex}(\mc{F})$. We will show that $f,\alpha_1,\alpha_2$ satisfy the theorem. By the choice of $f$, we have $\psi(x) = \min t(K_2; W)$, where the minimum is over all $W \in \mathrm{Lex}(\mc{F})$ with $t(P_4;W)=x$.

Consider $W= \otimes_{i=1}^\infty J_i$, where for every $i\ge 1$ either $J_i = (F_1,\mu)$ or $J_i=(F_2,\mu)$. Let $T = \{i \in \mathbb{N} \ : \ J_i = (F_1,\mu)\}$. Let  $\lambda = \sum_{i \in T} \frac{1}{2^{i}}$, and let $\gamma=\mathbf{m}_2(\mu) > 1/2$.
By Lemma~\ref{lem:density} we have
\begin{equation}\label{e:P4}
t(P_4; W) =\sum_{j=1}^\infty \frac{1}{2^{j-1}} t(P_4;J_j)=\sum_{i \in T} \frac{1}{2^{i-1}} \alpha_1 + \sum_{i \in \bb{N} - T} \frac{1}{2^{i-1}} \alpha_2 = \lambda \alpha_1 + (1-\lambda) \alpha_2,
\end{equation}
and, similarly,
$$t(K_2;W) = \sum_{j=1}^\infty\gamma^{j-1} t(K_2;J_j)= \beta_2 + (\beta_1-\beta_2)(1-\gamma)\sum_{i\in T}\gamma^{j-1}.$$
The identity (\ref{e:P4}) shows that when $\lambda$ is not a dyadic rational, and $x = \lambda \alpha_1 + (1-\lambda) \alpha_2$ a convex combination of $\alpha_1$ and $\alpha_2$, then there
exists a unique $W \in \mathrm{Lex}(\mc{F})$ that satisfies $t(P_4,W)=x$, and if $\lambda$ is dyadic then there are two such $W \in \mathrm{Lex}(\mc{F})$  corresponding to the two different binary expansions of $\lambda$.    Out of the two choices of   $(\lambda_i)_{i \in \mathbb{N}}$,  the one that corresponds to a finite $T$ minimizes $t(K_2;W)$.

Now it can be easily seen now that for every $x \in [\alpha_1,\alpha_2]$,
$$\limsup_{h \to 0} \frac{\psi(x+h)-\psi(x)}{h}= +\infty,$$
 finishing the proof of the theorem.
\end{proof}

\section{Concluding remarks}\label{sec:conc}

{\bf Differentiability of the boundary of $T(\mc{F})$.}
In the main theorem of this article, Theorem~\ref{thm:main}, we showed that there exists a quantum graph $f$ such that the set of points $(t(K_2;W),t(P_4;W))$ for all graphons $W$ that satisfy  $t(f;W)=0$ is the union of a curve that has a nowhere differentiable part and a countable set of points. The region defined by $f$ is of measure $0$, however it is not difficult to slightly modify $f$ to obtain a region of positive measure whose boundary is nowhere differentiable in some parts.

Let $f= \sum_{F \in \mc{F}} \lambda_F F$ for  be the quantum graph in Theorem~\ref{thm:main}. The theorem says that the projection of the intersection of $T(\mc{F} \cup \{K_2,P_4\})$ with the hyperplane defined by $\sum_{F \in \mc{F}} \lambda_F x_F=0$ to the coordinates $(x_{K_2},x_{P_4})$ is a curve that is nowhere differentiable  in certain parts plus a countable set of points. Unfortunately this does not imply anything about differentiability of the boundary of $T(\mc{F} \cup \{K_2,P_4\})$. The main question still remain open: For a finite collection $\mc{F}$ of graphs, is the boundary of the region $T(\mc{F})$ almost everywhere differentiable?

\vskip 5pt

\noindent
{\bf Finite forcibility of lexicographic powers.}
Lemma~\ref{lem:key} implies that  $\otimes^{\infty}F$ is finitely forcible for every stringent graph $F$, giving a partial answer to a question of Lov\'{a}sz and Szegedy~\cite{MR2802882} mentioned in Section~\ref{s:forcibility}. We believe that one can extend the argument to show that  $\otimes^{\infty}F$ is finitely forcible for every prime $F$. Dealing with homogeneous sets, however, presents a major technical issue.

\vskip 10pt
\noindent {\bf Acknowledgement.} The second author thanks Liana Yepremyan for discussions related to the subject of this paper.


\bibliographystyle{amsalpha}
\bibliography{diff}
\end{document}